\documentclass[]{amsart}
\usepackage{color}
\usepackage{amsmath}
\usepackage{amsthm}
\usepackage{amssymb}
\usepackage{graphicx}
\usepackage{amsfonts}
\usepackage{mathrsfs}

\theoremstyle{plain}
\newtheorem{theorem}{Theorem}[section]

\theoremstyle{definition}

\newtheorem{remark}[theorem]{Remark}
\newtheorem*{remark*}{Remark}

\newcommand{\ee}{\mathrm{e}}
\newcommand{\ii}{\mathrm{i}}

\def\dbar{\overline{\partial}}
\def\pmtwo#1#2#3#4{\left( \begin{array}{cc}#1&#2\\#3&#4\end{array}\right)}

\begin{document}

\title{Spectral approach to the scattering map  for the semi-classical
defocusing Davey-Stewartson II equation}

\author[C.~Klein]{Christian Klein}
\address[C.~Klein]
{Institut de Math\'ematiques de Bourgogne
9 avenue Alain Savary, BP 47870, 21078 Dijon Cedex}
\email{christian.klein@u-bourgogne.fr}

\author[K.~McLaughlin]{Ken McLaughlin}

\address[K.~McLaughlin]
{Department of Mathematics,
1874 Campus Delivery,
Fort Collins, CO 80523-1874}
\email{kenmcl@rams.colostate.edu}

\author[N.~Stoilov]{Nikola Stoilov}

\address[N.~Stoilov]{Institut de Math\'ematiques de Bourgogne, UMR 5584\\
                Universit\'e de Bourgogne-Franche-Comt\'e, 9 avenue Alain Savary, 21078 Dijon
                Cedex, France}
    \email{Nikola.Stoilov@u-bourgogne.fr}

\begin{abstract}
    The inverse scattering approach for the defocusing 
    Davey-Stewartson II equation is given by a system of D-bar 
    equations. 
We present a numerical approach to semi-classical D-bar problems for 
real analytic rapidly decreasing potentials. We treat the D-bar problem as a complex linear second order integral equation which is solved with discrete Fourier transforms complemented by a 
regularization of the singular parts by explicit analytic computation. The 
resulting algebraic equation is solved either by fixed point iterations or GMRES. 
Several examples for small values of the semi-classical parameter in 
the system are discussed. 
\end{abstract}

\date{\today}

\subjclass[2000]{}

\thanks{}
\maketitle

\section{Introduction}
We present an efficient numerical approach to the direct and inverse scattering 
maps for the defocusing Davey-Stewartson (DS) II equation. The latter appears in many 
applications in nonlinear optics and hydrodynamics and can be written 
in the form
\begin{equation}
\begin{array}{l}
\label{eq:DSIIa}
i \epsilon q_{t} + \frac{\epsilon^{2}}{2} \left( q_{xx} - q_{yy} \right) = - |q|^{2} q - \varphi q, \\
~\\
 \varphi_{xx} + \varphi_{yy} = -2 \left( |q|^{2} \right)_{xx}.
\end{array}
\end{equation}
This equation is an integrable 2d generalisation of the nonlinear Schr\"odinger (NLS) equation \cite{AH}. Here $\epsilon\ll 1$ 
is a small parameter similar to $\hbar$ in the Schr\"odinger 
equation, and the limit $\epsilon\to 0$ is accordingly called the semiclassical limit. Since the DS II equation 
is purely dispersive, in the semiclassical limit the solutions show  
zones of rapid modulated 
oscillations called \emph{dispersive shock waves} (DSW), see for instance 
\cite{KR,KR2}. 

The study of DSWs for generic dispersive partial differential 
equations (PDEs) is limited to multi-scales analysis for sufficiently 
small initial data. For completely integrable PDEs as the NLS 
equation, a more general description of DSWs can be given, see for 
instance \cite{BK,BT,DGK,JLM} and references therein. Since the DS II 
equation is completely integrable, see \cite{AH,Sch}, an  
approach similar to the treatment of NLS could be applied. A first step in this direction has 
been taken in \cite{AKMM}.

The inverse scattering map for the defocusing DS II equation is given 
by the following elliptic system.
\begin{eqnarray}
\label{eq:specprob2}
\epsilon\pmtwo{\dbar}{0}{0}{\partial} \psi = 
\frac{1}{2}\pmtwo{0}{q}{\overline{q}}{0} \psi,\quad \epsilon\ll1 .
\end{eqnarray}
The operators $\partial $ and $\overline{\partial}$ are defined via
\begin{eqnarray*}
\partial = \frac{1}{2} \left(
\frac{\partial}{\partial x} - i \frac{\partial}{\partial y}
\right), \ \ \ \ \ 
\overline{\partial} = \frac{1}{2} \left(
\frac{\partial}{\partial x} + i \frac{\partial}{\partial y}
\right).
\end{eqnarray*}
We note here that this system appears in a number of different applications, 
for instance in the context of further integrable partial 
differential equations \cite{BC}, of 2D orthogonal polynomials, of 
Normal Matrix Models in Random Matrix Theory, see e.g.\ \cite{KM}, 
and of electrical impedance tomography (EIT) \cite{CIN,MuSi2012,Uhl}. 
The D-bar system appearing in the inverse scattering approach to the 
defocusing Davey-Stewartson II equation
 has the most general form (see \cite{Perry2012} and 
references therein) and is hence 
the subject of this work. 

In all applications the quest is to recover a vector-valued solution $\psi = \psi(z,k) = \left( \begin{array}{c}
\psi_{1} \\
\psi_{2} \\ \end{array}
\right)$ with the following asymptotic behaviour as $|z| \to \infty$:
\begin{eqnarray*}
&&\lim_{|z| \to \infty} \psi_{1} e^{-kz/\epsilon} = 1, \\
&& \lim_{|z| \to \infty} \psi_{2} e^{-\overline{k} \overline{z} / \epsilon} = 0 \ ,
\end{eqnarray*}
writing $k=k_{1} + i k_{2} $ for $(k_{1},k_{2})\in\mathbb{R}^2$, a complex-valued parameter playing the role of a spectral variable.
The quantity $\psi$ is referred to as a \emph{complex geometric 
optics} (CGO) solution.  Information relevant to the inverse problem at hand is obtained from the CGO solution by extracting the value at either $0$ or $\infty$, both of which are delicate limits.  

For the DSII equation, the {\it reflection coefficient}, $r=r^{\epsilon}(k)$, is encoded in the sub-leading term in the asymptotic expansion of $\psi$ as $z \to \infty$, via
\begin{eqnarray}\label{reflec}
\psi_{2} e^{ - \overline{k} \overline{z}} = \frac{\bar{r^{\epsilon}}(k)}{2\overline{z}} + \mathcal{O} \left( \frac{1}{|z|^{2}} \right) \ .
\end{eqnarray}
In fact, the notation $r(k$) does not imply that the function is 
holomorphic (and the same holds for any function of a complex 
variable in this paper), and what we have is a transformation from the potential 
$q(x,y)$ (again a function of two real variables) to a function 
$r(k_{1}, k_{2})$, which extends to Lipschitz continuous and 
invertible mapping on the function space $L^{2}\left( \mathbb{C} 
\right)$ (see \cite{Perry2012,NRT}, and the references contained therein).  

Now, if $q = q(x,y,t,\epsilon)$ evolves according to the DSII equation (\ref{eq:DSIIa}), then the reflection coefficient evolves according to
\begin{eqnarray*}
r=r^{\epsilon}(k,t) = r^{\epsilon}(k,0) e^{ \frac{-it}{4 \epsilon}  \left( k^{2} + \overline{k}^{2} \right) } \ .
\end{eqnarray*}
In this sense, the mapping from $q$ to $r$ linearises the DSII flow.  
More amazingly, the inverse problem of reconstructing the potential $q(x,y,t,\epsilon)$ from the reflection coefficient $r(k,t)$ is also a D-bar problem, only in the complex variable $k$.  Indeed, setting 
\begin{equation}
\Phi_1=\Phi^\epsilon_1(k;z,t):=
\ee^{-kz/\epsilon}\psi_1\quad\text{and}\quad
\Phi_2=\Phi^\epsilon_2(k;z,t):=
\ee^{-\overline{k} \overline{z}/\epsilon}\psi_2 \ , 
\label{Phi}
\end{equation}
it turns out that one has, for each $z \in \mathbb{C}$, 
\begin{equation}
\epsilon\overline{\partial}_k \Phi_1=\tfrac{1}{2}e^{(\overline{k} 
\overline{z} - k z )/ \epsilon} \overline{r^\epsilon(k,t)}\Phi_2, \ \ \ \ \ \ \ \ 
\epsilon \partial_k \Phi_2 =\tfrac{1}{2}e^{-(\overline{k} \overline{z} - k z )/ \epsilon}r^\epsilon(k,t)\Phi_1
\label{eq:dbar-k}
\end{equation}
where,
\begin{equation*}
\overline{\partial}_k:=\frac{1}{2}\left(\frac{\partial}{\partial k_{1}}+\ii\frac{\partial}{\partial k_{2}}\right), \ \ \ \ \ \partial_{k}:=\frac{1}{2}\left(\frac{\partial}{\partial k_{1}}-\ii\frac{\partial}{\partial k_{2}}\right),
\end{equation*}
and the asymptotic conditions
\begin{equation}
\lim_{|k|\to\infty}\Phi^\epsilon_1(k;z,t)=1\quad\text{and}\quad
\lim_{|k|\to\infty}\Phi^\epsilon_2(k;z,t)=0.
\label{eq:Phi-asymp}
\end{equation}
The functions $\Phi_{1}$ and $\Phi_{2}$, being uniquely determined by the above elliptic system (\ref{eq:dbar-k}) and boundary conditions (\ref{eq:Phi-asymp}), yield the potential $q(x,y,t,\epsilon)$ through the asymptotic behavior as $|k| \to \infty$:
\begin{eqnarray*}
\Phi_{2} =  \frac{\overline{q(x,y,t,\epsilon)}}{2 k } + \mathcal{O} \left( |k|^{-2}\right) \ .
\end{eqnarray*}

In this paper we present a numerical method that is spectrally 
accurate to solve the D-bar equation (\ref{eq:specprob2}) with the desired asymptotic conditions.
The potential $q = q(x,y)$ is assumed to be in the Schwartz class 
$\mathcal{S}(\mathbb{R}^{2})$ of 
rapidly decreasing smooth functions.   

Knudsen, Mueller and Siltanen 
\cite{KnMuSi2004} developped a numerical approach to solving D-bar problems of the form
\begin{eqnarray*}
\overline{\partial} M = q(x,y) \overline{M}, \ \ \ \ \ M = 1 + 
\mathcal{O} \left(  \frac{1}{|z|}\right) \ \mbox{ as } |z| \to \infty \ . 
\end{eqnarray*}
They use Fourier techniques for these equations in an integral 
representation (see below). Their method is of first order since the singular 
integrand was regularized by putting it equal to zero where it 
diverges. In \cite{KM} the authors present an enhanced version of 
this approach having spectral accuracy for 
potentials $q(x,y)$ in the Schwartz class. They use this numerical method to compute the reflection coefficient $r(k)$, as well as the inverse problem of computing the potential $q(x,y,t)$, in the case that $\epsilon = 1$.  (The numerical method is summarized in Section~\ref{int_appr} below.)  Subsequent to that, in \cite{AKMM}, the asymptotic behavior of the DSII equation (\ref{eq:DSIIa}) with $\epsilon \to 0$ is considered, along with a number of different numerical methods aimed at elucidating the challenges in the asymptotic analysis.

In this work we are interested in the numerical construction of CGO 
solutions to  equation
(\ref{eq:specprob2}) for potentials $q$ in the Schwartz class 
$\mathcal{S}(\mathbb{R}^{2})$ of 
rapidly decreasing smooth functions, with $\epsilon$ tending to zero. The paper is organised as follows: In Section~\ref{int_appr} we reformulate the D-bar problem as an integral equation and discuss the basic ingredients of our numerical implementation. In Section~\ref{secfix} we solve the discretised integral equation of the preceding section by a fixed point iteration, whereas in Section~\ref{secgmres} we do the same by using the GMRES method. Finally, we draw our conclusions in Section~\ref{secOutlook}.     

\section{Integral equation and numerical approaches}\label{int_appr} 
In this section we present a reformulation of system 
(\ref{eq:specprob2}), which is suited for an efficient numerical 
treatment and discuss the employed numerical approach. 

\subsection{Integral equations}

The CGO solutions to (\ref{eq:specprob2}) have an essential 
singularity at infinity which is numerically problematic for obvious 
reasons. The quantities $\Phi_{1}$ and $\Phi_{2}$ defined in (\ref{Phi}), with asymptotic normalization (\ref{eq:Phi-asymp}), are well suited for numerical simulations.
As stated in the introduction, in terms of these functions, system (\ref{eq:specprob2}) takes the form
\begin{equation}
\begin{array}{l}
    \epsilon\bar{\partial}\Phi_{1}=\frac{1}{2}q\mathrm{e}^{(\bar{k}\bar{z}-kz)/\epsilon}\Phi_{2},\\
    ~\\
    \epsilon\partial\Phi_{2}=\frac{1}{2}\bar{q}\mathrm{e}^{(kz-\bar{k}\bar{z})/\epsilon}\Phi_{1}.
    \label{dbarphi}
\end{array}
\end{equation}
\begin{remark}
    System (\ref{dbarphi}) is related to the corresponding system 
    with $\epsilon=1$ via the transformation $q\mapsto q/\epsilon$ and 
    $k\mapsto k/\epsilon$. Since we vary $\epsilon$ in the examples, 
    we will concentrate on potentials with $||q||_{\infty}\sim 1$ 
    (smaller values of $\epsilon$ imply larger values of $|q|$ in the 
    system with $\epsilon=1$) and 
    discuss only the cases $k=0$ and $k=1$. 
\end{remark}

We define $\xi=\xi_{1}+i\xi_{2}$ where $\xi_{1}$ and $\xi_{2}$ are 
the dual Fourier variable to $x$ and $y$ respectively. The Fourier 
transform of a function $\Phi$ is defined as 
\begin{align*}
    \hat{\Phi} & = \mathcal{F}\Phi:=\frac{1}{2\pi}\int_{\mathbb{R}^{2}}^{}\Phi 
    e^{-i(\xi \bar{z}+\bar{\xi}z)/2}\Phi dx dy,
    \\
    \Phi & =\mathcal{F}^{-1}\Phi=\frac{1}{2\pi}\int_{\mathbb{R}^{2}}^{} e^{i(\xi 
    \bar{z}+\bar{\xi}z)/2}\hat{\Phi} d\xi_{1}d\xi_{2}.
\end{align*}
This implies formally
$$\mathcal{F}(\bar{\partial}\Phi_{1})=\frac{i}{2}\xi 
\hat{\Phi}_{1},\quad \mathcal{F}(\partial\Phi_{2})=\frac{i}{2}\bar{\xi} 
\hat{\Phi}_{2}.$$
Thus we get for the second relation in (\ref{dbarphi})
\begin{equation}
    \Phi_{2}=\mathcal{F}^{-1}\left[\frac{1}{i\bar{\xi}\epsilon}\mathcal{F}(\bar{q}\Phi_{1})\circ
    (\xi-2i\bar{k}/\epsilon)
    \right],
    \label{Phi21}
\end{equation}
which is equivalent to 
\begin{equation}
    \Phi_{2}=e^{(kz-\bar{k}\bar{z})/\epsilon}\mathcal{F}^{-1}\left[\frac{1}{i\epsilon(\bar{\xi}-2ik/\epsilon)}\mathcal{F}(\bar{q}\Phi_{1}) \right].   
    \label{Phi2osc}
\end{equation}
    Therefore $\Phi_{2}$ is given by a singular Fourier integral times an 
    oscillatory term.    
An analogous formula can be obtained for $\Phi_{1}$,
\begin{equation}
    \Phi_{1}=\mathcal{F}^{-1}\left[\frac{1}{i\xi\epsilon}\mathcal{F}(q\Phi_{2})\circ
    (\xi-2i\bar{k}/\epsilon)
    \right].
    \label{Phi11}    
\end{equation}¥
Replacing in (\ref{Phi11})  
$\Phi_{2}$ via (\ref{Phi2osc}), one gets an integral equation for 
$\Phi_{1}$,
\begin{equation}
    \Phi_{1}=-\mathcal{F}^{-1}\left\{\frac{1}{\epsilon\xi}\mathcal{F}\left[q\mathcal{F}^{-1}\left(
    \frac{\mathcal{F}(\bar{q}\Phi_{1})}{\epsilon(\bar{\xi}-2ik/\epsilon)}\right)\right]\right\}
    \label{Phi1},
\end{equation}
i.e., a singular integral equation for $\Phi_{1}$ which, once discretised, can be solved 
iteratively either by GMRES or by a fixed point iteration. Remarkably, 
this equation does not contain oscillatory terms, and it is complex 
linear in $\Phi_{1}$. After obtaining a solution of 
this equation, the function $\Phi_{2}$ follows from (\ref{Phi2osc}) 
via a singular integral. 

\begin{remark}
    System (\ref{eq:specprob2}) can be written as a pure D-bar system 
    by taking the complex conjugate of the second equation and by 
    working with $\bar{\psi}_{2}$ instead of $\psi_{2}$. 
    Diagonalising the matrix on the resulting right hand side yields 
    equations of the form ($\mu_{\pm}:=\psi_{1}\pm\bar{\psi}_{2}$)
    $$\bar{\partial}\mu_{\pm}=\pm\frac{q}{2}\bar{\mu}_{\pm},$$
    see for instance \cite{KM}. These equations are not complex 
    linear and thus have to be split into real and imaginary part in 
    order to apply the Krylov techniques \cite{gmres} we will use in 
    section \ref{secgmres}. This essentially doubles memory 
    requirements which makes the approach of the present paper much 
    more efficient in this respect. 
\end{remark}

However, equation (\ref{Phi1}) is not ideal for a numerical treatment 
with discrete Fourier transforms, since the functions $\Phi_{1,2}$ are not in the Schwarz class.  Indeed, the singular Fourier symbols 
appearing in the equations (\ref{Phi1}) and (\ref{Phi2osc}) cause these functions to decay too slowly at $\infty$.
But the function $S$  defined by (essentially the Fourier transform of the D-bar 
derivative of $\Phi_{1}$)
\begin{equation}
    S := \xi \hat{\Phi}_{1},\quad 
    \Phi_{1}=\mathcal{F}^{-1}\left(\frac{S}{\xi}\right) + 1
    \label{S},
\end{equation}
is in this class. Thus we will consider the equation following from 
(\ref{Phi1}),
\begin{equation}
    S=-\mathcal{F}\left[q\mathcal{F}^{-1}\left(
    \frac{1}{\epsilon(\bar{\xi}-2ik/\epsilon)}\mathcal{F}\left\{\frac{\bar{q}}{\epsilon}\left(
    \mathcal{F}^{-1}\left(\frac{S}{\xi}\right) + 1\right)\right\}\right)\right]
    \label{Seq}.
\end{equation}

The reflection coefficient (\ref{reflec}) was discussed in detail in 
\cite{KM}, here the focus is on the solutions to the D-bar system 
(\ref{eq:specprob2}). Nevertheless, we note here that it can be obtained from a 
given solution $\Phi_{1}$ in a straight forward way at essentially no 
additional computational cost. Note that the solution for $\Phi_{2}$ 
in (\ref{dbarphi}) can be written in terms of a solid Cauchy transform as 
\begin{equation*}
    \Phi_{2}=\frac{1}{\pi}    
    \int_{\mathbb{R}^{2}}^{}\frac{e^{(kz'-\bar{k}\bar{z}')/\epsilon}\bar{q}\Phi_{1}}{\epsilon(\bar{z}-\bar{z}')}
    dx'dy'
\end{equation*}
Thus we get with (\ref{reflec}) by computing $\lim_{|z|\to\infty }
\bar{z}\Phi_{2}$
\begin{equation}
    \bar{r}(k)=\frac{2}{\epsilon\pi}    
    \int_{\mathbb{R}^{2}}^{}e^{(kz-\bar{k}\bar{z})/\epsilon}\bar{q}\Phi_{1}
    dxdy
    \label{reflec3}.
\end{equation}
This, up to a multiplicative factor,
corresponds  to the Fourier transform of 
$e^{kz-\bar{k}\bar{z}}\bar{q}\Phi_{1}$ for $\xi=0$. 
Since we will work on a Fourier grid as detailed below, this integral 
can be  computed simply as the mean value of the integrand on the grid thus providing a spectral method, see the 
discussion in \cite{trefethen}. For given $\Phi_{1}$, one just 
has to multiply the function with 
$e^{kz-\bar{k}\bar{z}}\bar{q}$ and compute the mean value. 
For a potential in the Schwartz class, this approach is 
efficient even for large $k$ since $q$ is rapidly decreasing with 
$|z|$ thus delimiting the effects of the oscillatory term $e^{kz-\bar{k}\bar{z}}$. 

\subsection{Singular Fourier integrals}

\begin{figure}[htb!]
   \includegraphics[width=0.49\textwidth]{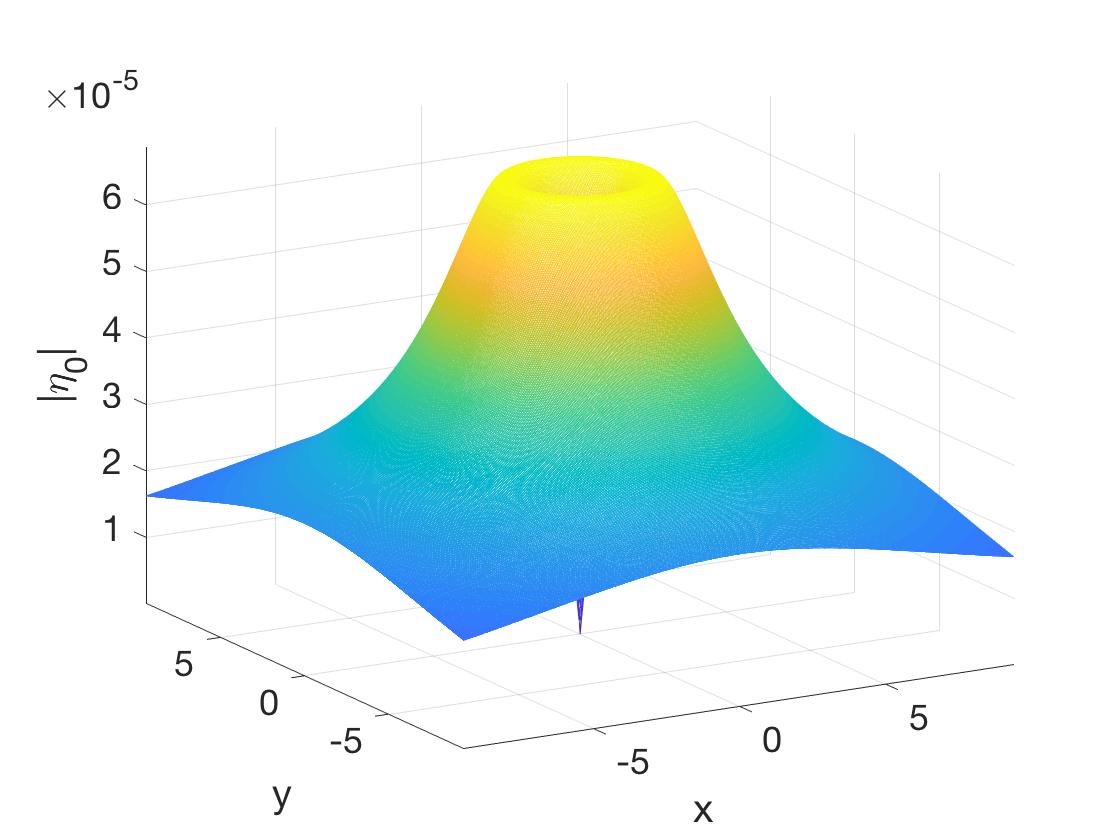}
   \includegraphics[width=0.49\textwidth]{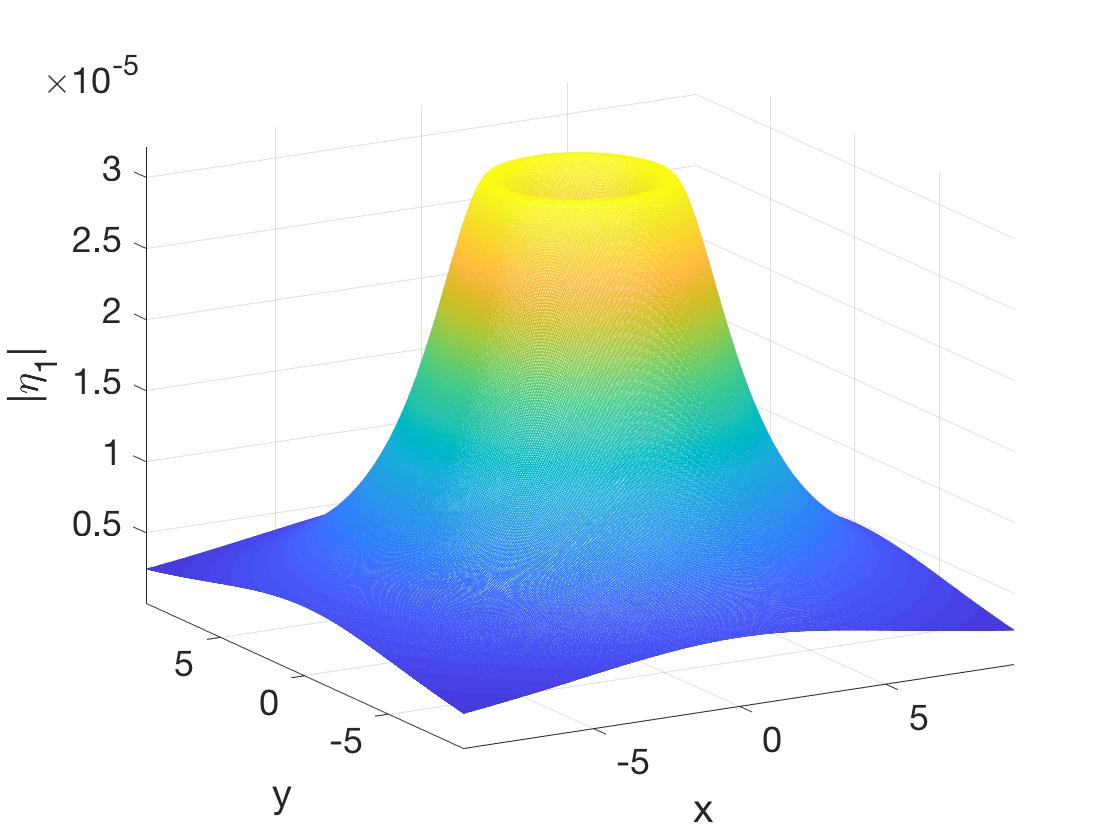}
\caption{The modulus of the functions $\eta_{0}$ (left) and 
$\eta_{1}$ (right) of (\ref{expxi2}).  }  
\label{eta}
\end{figure}

The task is thus to compute two singular integrals, the first being 
$\mathcal{F}^{-1}(S/\xi)$. As in \cite{KM} we observe that 
$2\mathcal{F}^{-1}\exp(-|\xi|^{2})=\exp(-|z|^{2}/4)$ and thus
\begin{equation}
    \mathcal{F}^{-1}\left(\frac{e^{-|\xi|^{2}}}{\xi}\right)=\frac{i}{z}\left(1-e^{-\frac{|z|^{2}}{4}}
    \right),
    \label{expxi}
\end{equation}
as well as 
\begin{equation}
    \mathcal{F}^{-1}\left(\frac{\bar{\xi}^{n}e^{-|\xi|^{2}}}{\xi}\right)=
    (-2i\partial)^{n}\frac{i}{z}\left(1-e^{-\frac{|z|^{2}}{4}}
    \right)=:\eta_{n},\quad n=0,1,\ldots
   \label{expxi2}
\end{equation}
We show the first two functions $\eta_{0}$, $\eta_{1}$ in 
Fig.~\ref{eta}. They vanish at the origin and have an annular 
structure. It can be seen that they are of order $10^{-5}$. The
$L^{\infty}$ norm of these functions decreases with $n$. 

The singular integral is then computed via
\begin{equation}
    \mathcal{F}^{-1}\left(\frac{S(\xi)}{\xi}\right)
    =     
    \mathcal{F}^{-1}\left(\frac{S(\xi)-e^{- |\xi|^{2}}\sum_{n=0}^{M}\partial_{\bar{\xi}}^{n}
    S(0)\bar{\xi}^{n}/n!}{\xi}\right)+ e^{- |\xi|^{2}} \sum_{n=0}^{M}\partial_{\bar{\xi}}^{n}
    S(0) (-2i\partial)^{n}\frac{i}{z}\left(1-e^{-\frac{|z|^{2}}{4}}
    \right).
    \label{xireg1}
\end{equation}
Here the integer $M$ is chosen  so that that the first term on the right hand side of the above equation is 
analytic within numerical precision, that is, its Fourier coefficients decrease to machine precision. 
The derivatives of $S$ are again computed via Fourier techniques,
\begin{equation}
    \partial_{\bar{\xi}}^{n}
    S(\xi)=\mathcal{F}\left[(-iz/2)^{n}\mathcal{F}^{-1}S\right]
    \label{xireg2}.
\end{equation}
Note that the derivatives are only needed for $\xi=0$. To compute the 
first term on the right-hand-side of (\ref{xireg1}) at $\xi=0$, we also need $\partial_{\xi}S(0)$ (according to L'Hopital's rule), which is computed as the derivatives in 
(\ref{xireg2}). 

The second singular integral $\mathcal{F}^{-1}[f(\xi)/(\bar{\xi}-2ik)]$ 
will be computed in essentially the same way. Note that we have 
\begin{equation*}
    \mathcal{F}^{-1}\left(\frac{e^{-|\xi|^{2}}}{\bar{\xi}}\right)=\frac{i}{\bar{z}}\left(1-e^{-\frac{|z|^{2}}{4}}
    \right),   
\end{equation*}
and a standard Fourier shifting calculation yields
\begin{equation*}
\mathcal{F}^{-1} \left( \frac{e^{ - | \xi + 2 i \overline{k} / \epsilon |^{2}}}{\overline{\xi} - 2 i k / \epsilon}\right) = e^{ - (kz - \overline{k} \overline{z} )/\epsilon} \frac{i}{\overline{z}} \left( 1- e^{-|z|^{2}/4} \right) \ ,
\end{equation*}
along with
\begin{equation*}
    \mathcal{F}^{-1}\left(\frac{ \left( \xi + 2 i \overline{k}/\epsilon\right) ^{n}e^{-|\xi + 2 i \overline{k}/\epsilon |^{2}}}{\bar{\xi} - 2 i k / \epsilon}\right)=e^{ - ( k z - \overline{k}\overline{z})/\epsilon}
    (-2i\bar{\partial})^{n}\frac{i}{\bar{z}}\left(1-e^{-\frac{|z|^{2}}{4}}
    \right).
\end{equation*}
We compute
\begin{align}
    \mathcal{F}^{-1}\left(\frac{f(\xi)}{\bar{\xi}-2ik/\epsilon}\right)
    &=     
    \mathcal{F}^{-1}\left(\frac{f(\xi)-
    e^{- |\xi + 2 i \overline{k} / \epsilon|^{2}}\sum_{n=0}^{M}\partial_{\xi}^{n}
    f(-2i\bar{k}/\epsilon)/n!(\xi+2i\bar{k}/\epsilon)^{n}}{\bar{\xi}-2ik/\epsilon}\right)\nonumber\\
    &+e^{(\bar{k}\bar{z}-kz)/\epsilon}
    \sum_{n=0}^{M}\partial_{\xi}^{n}
    f(-2i\bar{k}/\epsilon)/n! (-2i\bar{\partial})^{n}\frac{i}{\bar{z}}\left(1-e^{-\frac{|z|^{2}}{4}}
    \right),
    \label{xibarreg1}
\end{align}
where we have used the shifting argument in Fourier space from (\ref{Phi2osc}). 
The derivatives of $f$ are again computed via Fourier techniques,
\begin{equation}
    \partial_{\xi}^{n}
    f(\xi)=\mathcal{F}\left[(-i\bar{z}/2)^{n}\mathcal{F}^{-1}f\right]
    \label{xibarreg2}.
\end{equation}
The above regularization procedure is only applied if $2i\bar{k}/\epsilon$ is 
close to a $\xi$ grid point which means that the minimum of
$|\xi+2i\bar{k}/\epsilon|$ for $\xi$ a grid point is smaller than the minimal 
distance between two points of this grid. If this is the case, we denote 
this grid point by $\xi_{0}$ and apply the above formulae for 
$\xi_{0}$ in place of $-2i\bar{k}/\epsilon$. The value of the first term of 
the right hand side of (\ref{xibarreg1}) at this $\xi_{0}$ is then 
computed via $\partial_{\bar{\xi}}f(\xi_{0})$ where the derivative is 
obtained with the same method as in (\ref{xibarreg2}).

\subsection{Numerical implementation}
The standard numerical approach to approximate Fourier transforms 
as in (\ref{Seq}) is via discrete Fourier transforms which can be 
conveniently computed with a Fast Fourier transform (FFT). However, 
such methods will only be at their optimal efficiency when 
approximating periodic smooth functions. Functions in the Schwartz 
class are clearly not periodic, but when considered on a sufficiently 
large interval, they will have derivatives vanishing with the finite 
numerical precision at the boundaries of the computational interval. 
Thus the finite precision allows to approximate Schwartz class 
functions on such an interval as smooth periodic functions, and the 
coefficients of the FFT of the function will decrease exponentially 
as known for the standard Fourier transform of the function. Since also
this implies that  the numerical error decreases exponentially, 
such methods exhibit \emph{spectral convergence}. 

\begin{remark}\label{remschwartz}
    The problem of the D-bar system (\ref{dbarphi}) is that singular 
    Fourier multipliers appear which implies that the functions 
    $\Phi_{1,2}$ will not be in the Schwartz class, but have an algebraic 
    decrease in $z$ for $|z|\to\infty$ in contrast to their derivatives. 
    This is even a feature of crucial importance since the reflection 
    coefficient (\ref{reflec}) providing the scattering data in the context 
    of the DS equation is given as the coefficient of the $1/\bar{z}$ 
    term near infinity. This apparent incompatibility of the D-bar 
    solutions with an efficient setting for the application of FFT 
    techniques is addressed in this paper in the following way: 
    numerically only FFTs of smooth functions (within numerical 
    precision) are computed, all other terms will be computed by hand 
    as outlined in the previous subsection. 
\end{remark}

For the implementation we choose the computational domain
\begin{equation*}
    x\in L_{x}[-\pi,\pi], \quad y\in L_{y}[-\pi,\pi],
\end{equation*}
and the wave numbers in the FFT as 
\begin{equation*}
    \xi_{1}=(-N_{x}/2+1,-N_{x}/2+2,\ldots-N_{x}/2)/L_{x},\quad 
    \xi_{2}=(-N_{y}/2+1,-N_{y}/2+2,\ldots-N_{y}/2)/L_{y},
\end{equation*}
where $N_{x}$ and $N_{y}$ are the number of Fourier modes in the $x$ 
and $y$ direction respectively. As mentioned in remark 
\ref{remschwartz}, the values of $L_{x}$, $L_{y}$ and of $N_{x}$ and $N_{y}$ are 
chosen in a way that the Fourier transform of the studied Schwartz 
function decreases to machine precision as $\xi_{1}$ and $\xi_{2}$ 
approach the boundary of the computational domain.   Since the 
derivatives of $S$ in the regularization approach of the previous 
subsection are computed via (\ref{xireg2}) and (\ref{xibarreg2}), we 
want to ensure that also the inverse FFT of $S$ decreases to machine 
precision ($10^{-16}$ in our case, but in practice typically limited 
to $10^{-14}$ because of rounding errors) to limit also the numerical 
errors in the computation of $\mathcal{F}^{-1}S$ to this order. Thus 
we will show in the following examples that both the modulus of $S$ and 
of $\mathcal{F}^{-1}S$ decrease to machine precision for the chosen 
values of the parameters. 

The small size of the functions $\eta_{n}$, $n=0,1,2,\ldots$ of 
(\ref{expxi2}) implies that in practice only a small number $M$ of such 
functions is needed to compute the singular integrals to 
machine precision. For the examples studied in this paper, $M=4$ 
proved to be sufficient and is applied throughout all 
computations. 

\section{Fixed point iteration}\label{secfix}
The task to solve the D-bar system (\ref{dbarphi}) with the 
asymptotic conditions (\ref{eq:Phi-asymp}) has been reduced in the 
previous section to the 
solution of the singular integral equation (\ref{Seq}) and two 
singular quadratures (\ref{S}) and (\ref{Phi21}). In this section we 
will solve the integral equation with a fixed point iteration and study 
the convergence of the iteration for several examples. 

Thus we will solve equation (\ref{Seq}) iteratively in the form 
\begin{equation}
    S^{(n+1)}=-\mathcal{F}\left[q\mathcal{F}^{-1}\left(
    \frac{1}{\epsilon(\bar{\xi}-2ik/\epsilon)}\mathcal{F}\left\{\frac{\bar{q}}{\epsilon}\left(
    \mathcal{F}^{-1}\left(\frac{S^{(n)}}{\xi}\right) + 1\right)\right\}\right)\right]
    \label{Seqit},
\end{equation}
where $n=0,1,2,\ldots$ and $S^{(0)}=0$. The iteration is stopped once 
$||S^{(n+1)}-S^{(n)}||_{\infty}$ is smaller than some threshold which 
is typically taken to be $10^{-12}$. 

\begin{remark}
    The computational cost per iteration without regularization is 
    dominated by 4 2d FFTs. Per iteration there is a maximum of two 
    regularizations for the singular integrals if $2ik$ is close to a 
    point of the $\xi$ grid. Each of the regularizations requires an 
    additional 2d FFT to compute the derivatives of the integrand. But 
    since these derivatives are only computed at one single point, we 
    compute the inverse FFT of the  functions that need to be differentiated. 
    The derivatives are then sums of this function multiplied with 
    the respective powers of $z$ or $\bar{z}$ and an exponential factor. 
    The main computational cost is thus 6 2d FFTs per iteration. 
\end{remark}

As concrete examples we will study in this paper the potentials
\begin{align}
    q&=\exp(-x^{2}-y^{2})
    \label{gauss},\\
    q&=\exp(-x^{2}-3xy-5y^{2})
    \label{q2},    
\end{align}
i.e., a Gaussian and a rapidly decreasing potential which is not   
radially symmetric even asymptotically. 

\begin{figure}[htb!]
   \includegraphics[width=0.49\textwidth]{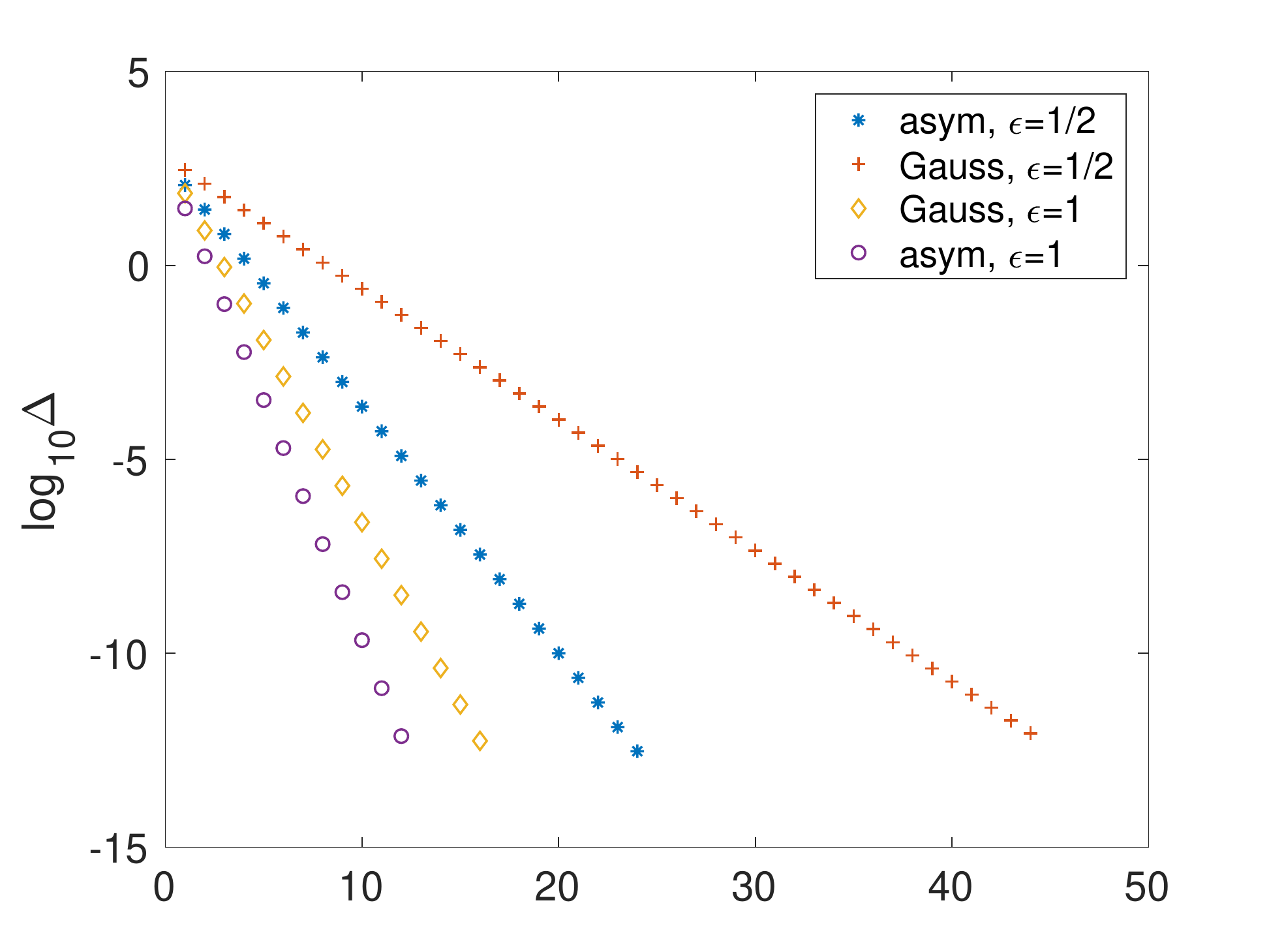}
   \includegraphics[width=0.49\textwidth]{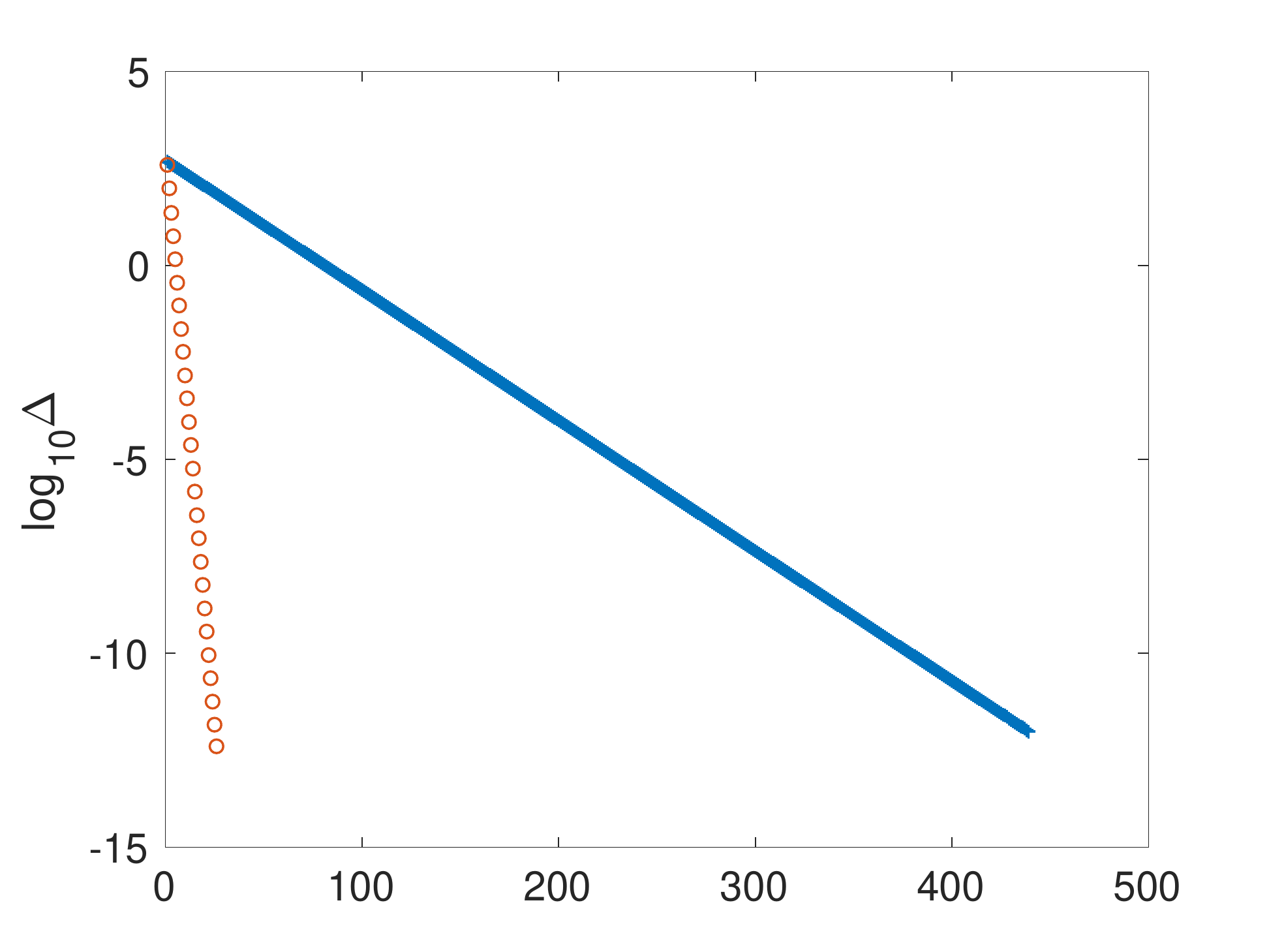}
\caption{The quantity $\Delta = ||S^{(n+1)}-S^{(n)}||_{\infty}$ in a 
logarithmic plot for the iterative solution of (\ref{Seqit}) for the 
potentials (\ref{gauss}) and (\ref{q2}) for $k=0$ and 
$\epsilon=1,1/2$ on the left; on the right the same quantity for the 
potential (\ref{q2}), $\epsilon=1/4$ and $k=0$ (`+') and $k=1$ (`o').  }  
\label{fix1}
\end{figure}

For the implementation we use $N_{x}=N_{y}=2^{8}$ Fourier modes and 
$L_{x}=L_{y}=3$. Putting $\Delta := ||S^{(n+1)}-S^{(n)}||_{\infty}$, 
one gets that the iteration converges linearly for $k=0$ and $\epsilon=1,1/2$ 
for the examples (\ref{gauss}), (\ref{q2}) as can be seen in 
Fig.~\ref{fix1} on the left. The  residuals computed for equation 
(\ref{Seq}) are $1.69*10^{-11}$, $6.71*10^{-12}$, $1.11*10^{-12}$, $8.39*10^{-13}$ as shown from left 
to right in the figure.  It can be noted that the convergence becomes 
slower the smaller $\epsilon$ is since this implies a potential $q$ 
of larger $L^{\infty}$ norm. Interestingly the convergence is slower 
for the symmetric (Gaussian) than for the non-symmetric potential. 
This can be seen on the right of Fig.~\ref{fix1} where the quantity 
$\Delta$ is shown for the potential (\ref{q2}) and $k=0,1$, 
$\epsilon=1/4$ (here we use $N_{x}=N_{y}=2^{9}$). Note that the iteration does not converge for the 
Gaussian potential for $k=0$ and $\epsilon=1/4$. In order to have a convergent 
scheme, the $L^{\infty}$ norm of the right hand side of 
(\ref{Seqit}), which is proportional to $1/\epsilon^{2}$, needs to 
be smaller than 1. On the other hand the iteration 
converges faster the larger $|k|$ as can be seen from the right 
figure of Fig.~\ref{fix1}. This is not surprising since a factor $k$ 
appears in the denominator of equation (\ref{Seqit}), see 
\cite{joh} for a more detailed discussion of the large $|k|$ 
limit. 

\begin{figure}[htb!]
   \includegraphics[width=0.49\textwidth]{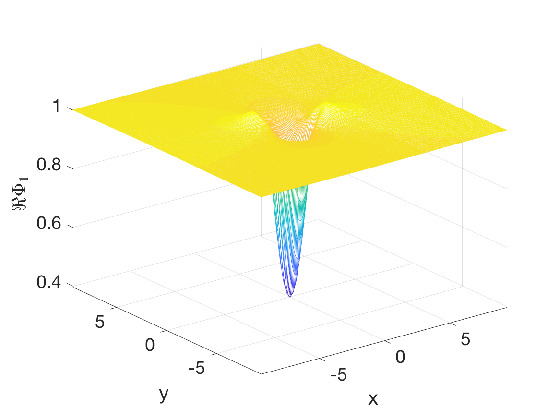}
   \includegraphics[width=0.49\textwidth]{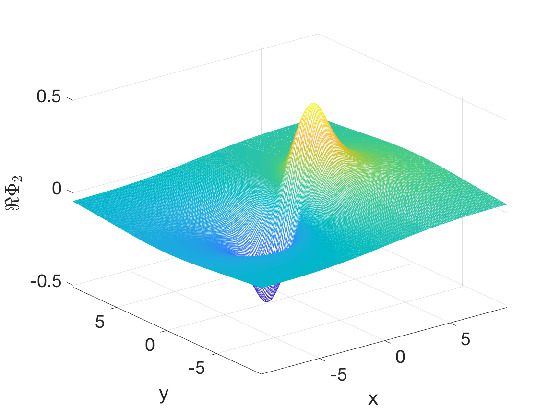}\\
   \includegraphics[width=0.49\textwidth]{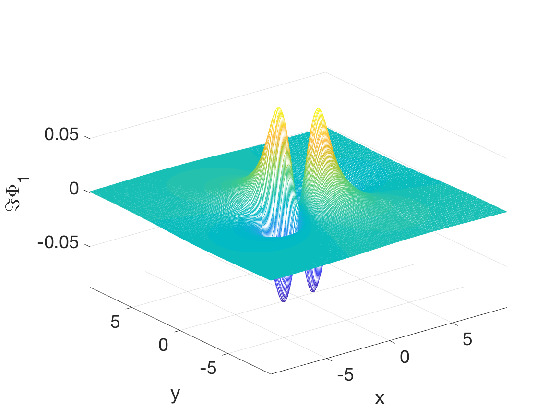}
   \includegraphics[width=0.49\textwidth]{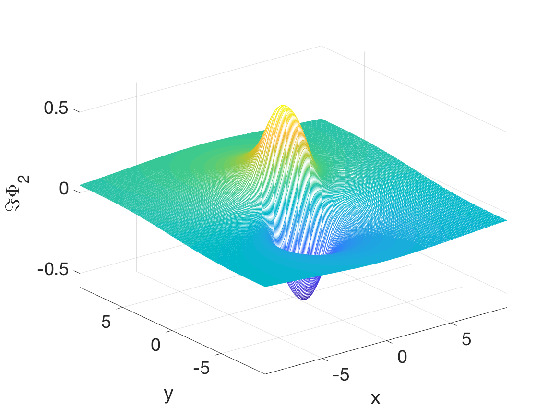}
\caption{Solutions to the system (\ref{dbarphi})  for the 
potential  (\ref{q2}) for $k=0$ and 
$\epsilon=1/4$, on the left $\Phi_{1}$, on the right $\Phi_{2}$, in 
the upper row the real parts, in the lower row the imaginary parts.  }  
\label{fix2}
\end{figure}

The solution corresponding to the right figure of Fig.~\ref{fix1} can 
be seen in Fig.~\ref{fix2}. Obviously the solutions tend only very 
slowly to the asymptotic values for $|z|\to\infty$. This comes again 
to show that it is only through the regularization approach we 
described in section \ref{int_appr} that Fourier methods can be used 
efficiently here. 

\begin{figure}[htb!]
   \includegraphics[width=0.49\textwidth]{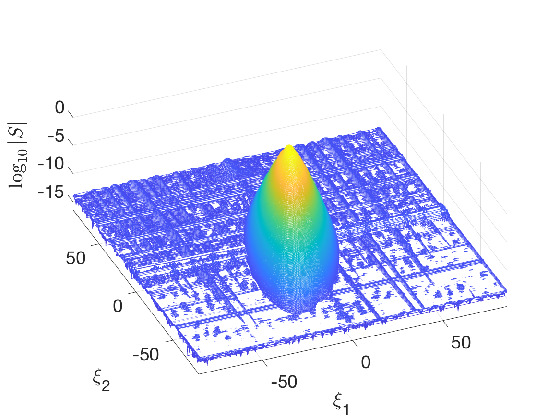}
   \includegraphics[width=0.49\textwidth]{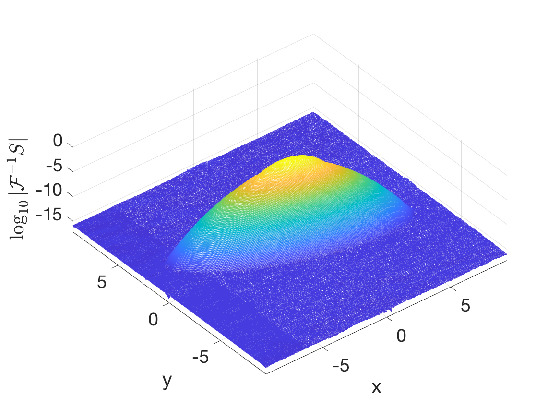}
\caption{Solution to equation (\ref{Seq})  for the 
potential  (\ref{q2}) for $k=0$ and 
$\epsilon=1/4$ on the left, and its inverse Fourier transform on the 
right. }  
\label{fix3}
\end{figure}

Since only the quantity $S$ is computed with FFT techniques, the 
numerical accuracy  can be controlled as mentioned in remark 
\ref{remschwartz} via the decrease of $|S|$ with $|\xi|$. A logarithmic 
plot of $|S|$ for the situation shown in Fig.~\ref{fix2} is 
presented on the left of Fig.~\ref{fix3}. It can be seen that $S$ 
decreases to machine precision. In order to control the numerical 
accuracy of the computation of the derivatives in (\ref{xireg2}) and 
(\ref{xibarreg2}) in the same way,  $\mathcal{F}^{-1}S$ also has to decrease to 
machine precision which is the case as shown on the right of 
Fig.~\ref{fix3}. 

\begin{figure}[htb!]
   \includegraphics[width=0.49\textwidth]{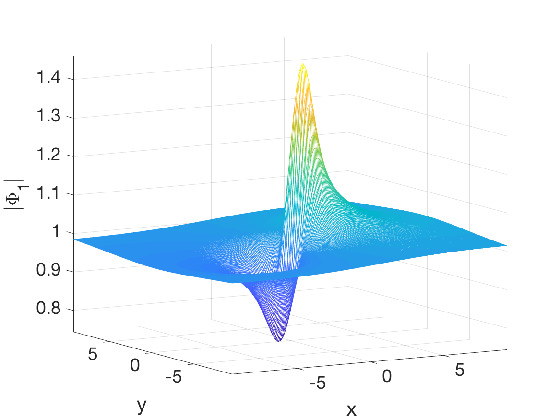}
   \includegraphics[width=0.49\textwidth]{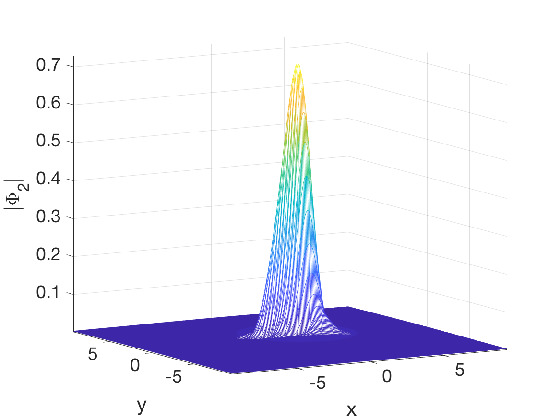}
\caption{Solutions to the system (\ref{dbarphi})  for the 
potential  (\ref{q2}) for $k=1$ and 
$\epsilon=1/4$, on the left $\Phi_{1}$, on the right $\Phi_{2}$.  }  
\label{fix4}
\end{figure}

The solution to the system (\ref{dbarphi}) for the potential 
(\ref{q2}) and $k=1$, $\epsilon=1/4$, for which the convergence was 
studied on the right of Fig.~\ref{fix1}, can be seen in 
Fig.~\ref{fix4}.

\section{Solution via GMRES}
\label{secgmres}
The fixed point iteration studied in the previous section converges 
without problems for potentials $q$ with an $L^{\infty}$ norm of 
order $\epsilon$. However this is no longer so for much larger $L^{\infty}$ norms, in which case the iteration becomes computationally expensive or diverges or both. Since we are interested 
as in \cite{AKMM} in the semi-classical limit $\epsilon\ll1$ (where 
$q$ is essentially replaced by $q/\epsilon$ in the equations with 
$\epsilon=1$), we 
present in this section an approach using GMRES. The fixed point 
iteration is mainly intended to provide a test of the accuracy of the 
GMRES approach to (\ref{Seq}) and to compare the convergence of the 
latter with the fixed point iteration.

The basic idea of using GMRES is that an equation of the form $Ax=b$, 
$x,b\in \mathbb{R}^{n}$ and $A$ an $n\times n$ matrix, is solved 
iteratively by approaching $A^{-1}b$ via linear combinations of $b$, 
$Ab,\ldots,A^{N}b$. The convergence of this approach is not 
guaranteed, and often a \emph{preconditioner} is needed, i.e.,  an $n\times n$ matrix $C$ such that 
the GMRES approach for $(CA)x=Cb$ converges. We do not use 
preconditioners here. The comparison with the fixed point iteration 
is used to highlight possible advantages of GMRES. 

An attractive feature of GMRES is that just the action of the matrix $A$ on a given 
vector is needed, not the matrix itself, thus allowing to proceed just storing $n$-dimensional vectors  instead of the 
$n\times n$ matrix $A$ thus reducing memory usage, a crucial advantage for the 
demanding problems studied here (note that the fixed point iterations 
has the same memory requirements as the GMRES approach). 
In our example, equation (\ref{Seq}) is to be solved after an FFT 
discretisation. We put
\begin{equation}
    b :=    \mathcal{F}\left[q\mathcal{F}^{-1}\left(
    \frac{1}{\bar{\xi}-2ik/\epsilon}\mathcal{F}\left\{\bar{q}\right\}\right)\right]
    \label{b}
\end{equation}
and
\begin{equation}
        AS:=\epsilon^{2}S+\mathcal{F}\left[q\mathcal{F}^{-1}\left(
    \frac{1}{\bar{\xi}-2ik/\epsilon}\mathcal{F}\left\{\bar{q}\left(
    \mathcal{F}^{-1}\left(\frac{S}{\xi}\right) \right)\right\}\right)\right]
    \label{A}.
\end{equation}
The matrix $b$ is written as a vector by putting 
the columns one after the other to form a vector of length 
$N_{x}N_{y}$, and in an analogous way for $AS$ (that is, column major ordering). 
Note that equation (\ref{Seq}) has been multiplied by $\epsilon^{2}$. 
This has no influence on the convergence of GMRES since the euclidean 
norm
of $Ax-b$ divided by the euclidean norm of $b$ controls the convergence. However, the 
absolute residuals  which are limited by machine precision can 
be chosen similar in this case to what could be imposed for the fixed point 
iteration, see the discussion below. 

\begin{figure}[htb!]
   \includegraphics[width=0.49\textwidth]{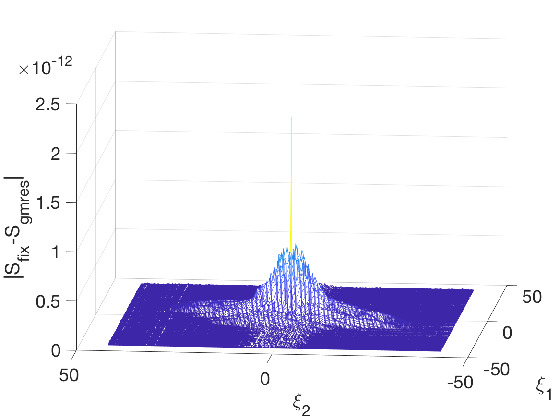}
   \includegraphics[width=0.49\textwidth]{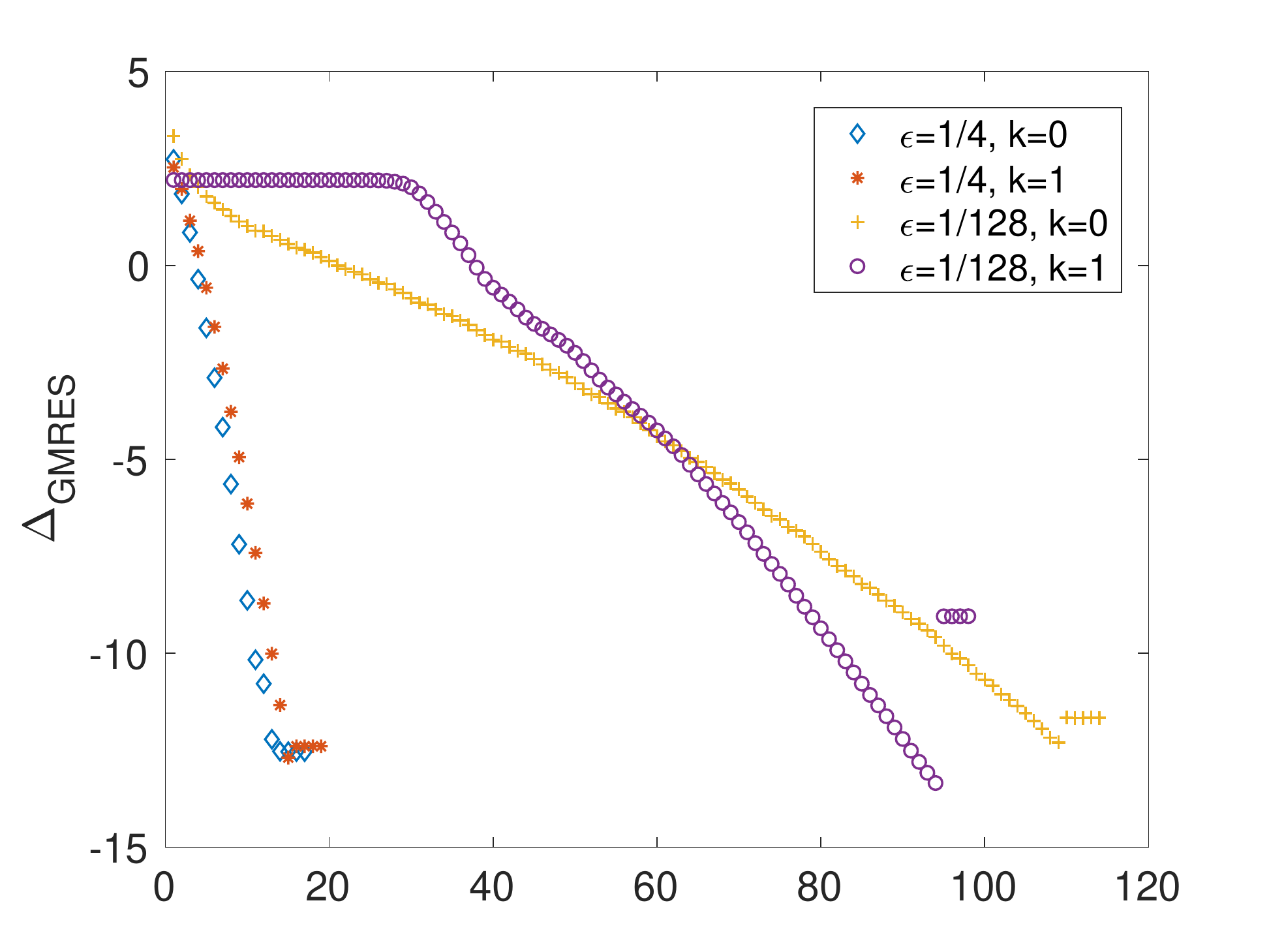}
\caption{Difference of the solutions to the system (\ref{dbarphi})  for the 
potential  (\ref{q2}) for $k=1$ and 
$\epsilon=1/4$ obtained with a fixed point iteration and GMRES on the 
left, and on the right the residuals $\Delta_{\mbox{GMRES}}$ for 
$\epsilon=1/4$, $\epsilon=1/128$ and $k=0,1$.  }  
\label{gmres1}
\end{figure}

We first study the situation of Fig.~\ref{fix4} with the same choice 
of the  parameters, but this time with GMRES. The iteration 
is stopped here once the relative residual drops below $10^{-14}$ 
which gives a residual of the same order as before. The difference between the 
solution with a fixed point iteration and the one with GMRES is shown 
in Fig.~\ref{gmres1} on the left. It is obviously largest for small 
$|\xi|$ (note that $S$ is in the Schwartz class), but overall of the 
order of $10^{-12}$ as expected. 

The right figure of 
Fig.~\ref{gmres1} shows the residuals $\Delta_{\mbox{GMRES}}$ in dependence of the 
number of iterates for various values of $\epsilon$ and $k$. It can 
be seen that for $\epsilon=1/4$, the approach reaches maximal 
accuracy for roughly 15 iterates which has to be compared to the 
corresponding figure on the right of Fig.~\ref{fix1} for the fixed 
point iteration. There for $\epsilon=1/4$ roughly 450 iterations were 
needed for $k=0$ and 40 for $k=1$. This shows that GMRES is much more 
efficient in this case. 

\begin{figure}[htb!]
   \includegraphics[width=0.49\textwidth]{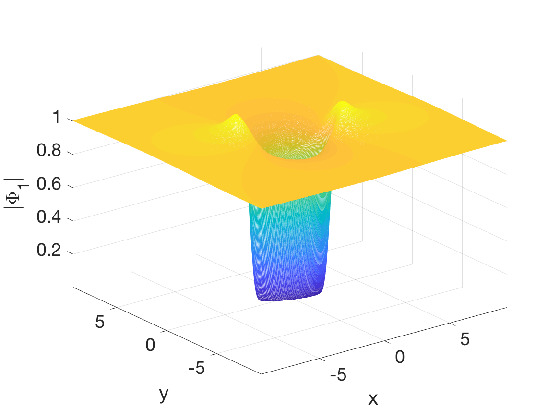}
   \includegraphics[width=0.49\textwidth]{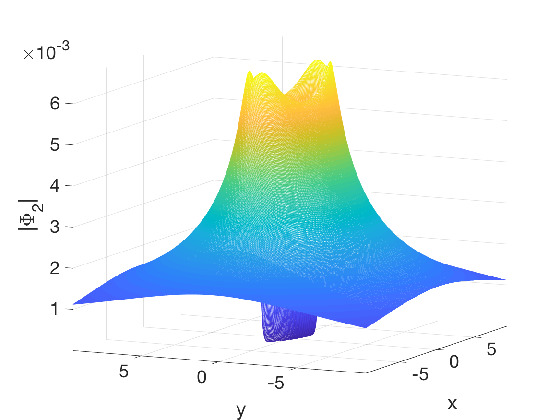}
\caption{Solutions to the(\ref{Seq})  for the 
potential  (\ref{q2})for $k=0$ and 
$\epsilon=1/128$ on the left $\Phi_{1}$, on the right $\Phi_{2}$.  }  
\label{gmres2}
\end{figure}

GMRES also allows to reach much smaller values of $\epsilon$ than 
accessible with the fixed point iteration. 
We use $L_{x}=L_{y}=3$ and $N_{x}=N_{y}=2^{10}$ Fourier modes to study 
the solution for the potential (\ref{q2}) for $\epsilon=1/128$ and 
$k=0$. The solutions can be seen in Fig.~\ref{gmres2}. Both solutions 
are clearly not radially sysmmetric, but are essentially constant in 
the vicinity of the origin. This is similar to the expected behavior 
for radially symmetric potentials discussed in \cite{AKMM}. GMRES 
converges after roughly 110 iterations as can be seen on the right of 
Fig.~\ref{gmres1} when a plateau reached. The relative residual is of 
the order of $10^{-15}$, the residual of $Ax-b$ with $A$, $b$ of
(\ref{A}), (\ref{b}) is of the order of $10^{-13}$. 

\begin{figure}[htb!]
   \includegraphics[width=0.49\textwidth]{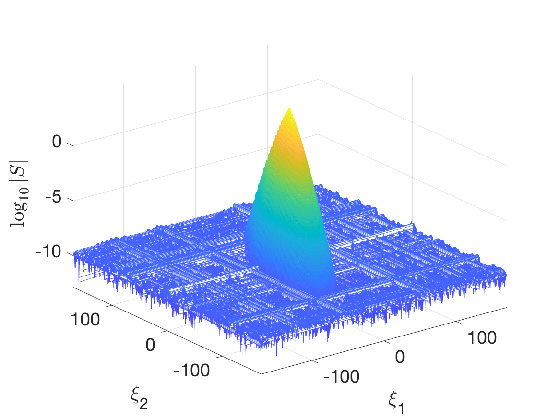}
   \includegraphics[width=0.49\textwidth]{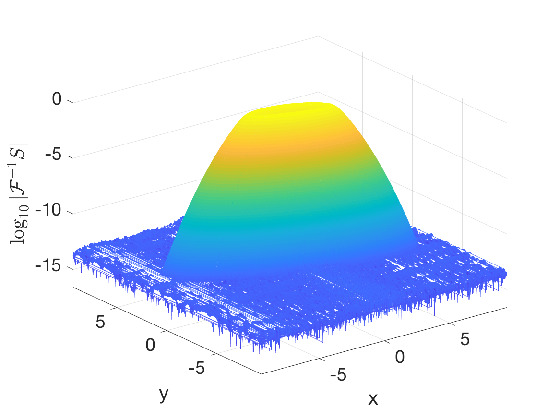}
\caption{Solution $S$  to the equation (\ref{Seq})  for the 
potential  (\ref{q2}) for $k=0$ and 
$\epsilon=1/128$ on the left and its inverse Fourier trasnform on the 
right.  }  
\label{gmres3}
\end{figure}

The parameters for the computation are chosen in a way that both $S$ 
and $\mathcal{F}^{-1}S$ decrease to the level of the rounding error 
as can be seen in Fig.~\ref{gmres3}. Note that the latter is of the order of 
$10^{-10}$ for the former since the maximum of $|S|$ is of the order 
of $10^{4}$ (in double precision it is in practice impossible to 
cover more than 14 orders of magnitude). 

\begin{figure}[htb!]
   \includegraphics[width=0.49\textwidth]{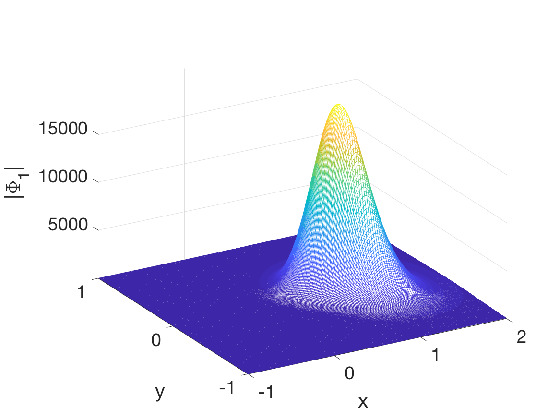}
   \includegraphics[width=0.49\textwidth]{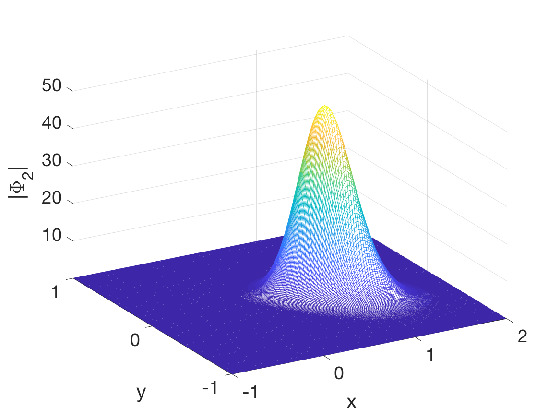}
\caption{Solutions to the system (\ref{dbarphi})  for the 
potential  (\ref{q2}) for $k=1$ and 
$\epsilon=1/128$; on the left $\Phi_{1}$, on the right $\Phi_{2}$.  }  
\label{gmres4}
\end{figure}

To study the solution for the potential (\ref{q2}) for 
$\epsilon=1/128$ and $k=1$, we use $N_{x}=2^{9}$ and $N_{y}=2^{10}$ 
Fourier modes and $L_{x}=2$, $L_{y}=1$. The solution can be seen in 
Fig.~\ref{gmres4}. Note that the solution $\Phi_{1}$ has a maximum of 
the order of $10^{4}$, but that the reflection coefficient being a 
rapidly decreasing function of the spectral parameter $|k|$ is of the 
order of $10^{-10}$ in this case. GMRES converges after 
roughly 90 iterations as can be seen on the right of 
Fig.~\ref{gmres1}. It stops since the last few iterations lead to a 
residual at a slightly higher plateau indicating that further 
iterations do not improve the accuracy. The residual of the equation is of 
the order of $10^{-12}$ when the iteration is stopped. 

\begin{figure}[htb!]
   \includegraphics[width=0.49\textwidth]{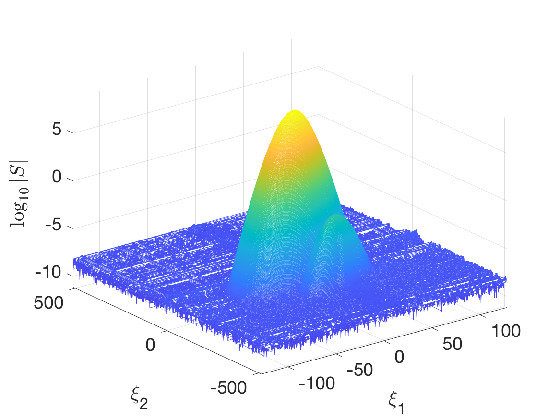}
   \includegraphics[width=0.49\textwidth]{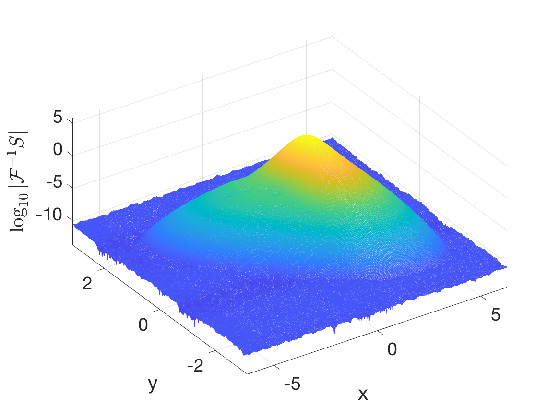}
\caption{Solution to equation (\ref{Seq})  for the 
potential  (\ref{q2}) for $k=1$ and 
$\epsilon=1/128$ on the left, and its inverse Fourier transform on the right.  }  
\label{gmres5}
\end{figure}

The numerical parameters are chosen in a way that $S$ and its inverse 
Fourier transform decrease to 
optimal precision. This can be seen on the left of Fig.~\ref{gmres5}. 
The function $S$ decreases to the order of the saturation level which 
is here at around $10^{-8}$ since the maximum of $|S|$ is of the 
order of $10^{6}$.  The inverse Fourier 
transform of $S$ is shown on right of Fig.~\ref{gmres5} to decrease 
to machine precision.

\section{Outlook}\label{secOutlook}
In this paper it was shown that CGO solutions to system 
(\ref{eq:specprob2}) for potentials in the Schwartz class 
can be efficiently constructed via the integral equation 
(\ref{Seq}) with FFT and iterative techniques. This can be done for a 
wide range of values of the spectral parameter $k$ and for the 
semiclassical parameter $\epsilon$. The limiting factor for small 
$\epsilon$ appears to be the conditioning of the matrix $A$ in 
(\ref{A}) which becomes worse the smaller $\epsilon$ is. We could 
reach values of $\epsilon=1/256$, but the achievable accuracy drops 
to the order of $10^{-5}$ in this case. 

This behavior is due to the singular character of the semiclassical 
solution which is discussed in \cite{AKMM}. There it is conjectured 
that the main contribution to the solution in this case to the CGO 
solutions is of the form $\exp(f/\epsilon)$ where $f$ solves an 
eikonal-type equation. It is not surprising that numerical approaches 
will ultimately fail to catch such a behavior for very small values 
of $\epsilon$. To address such questions, it appears best to use a 
\emph{hybrid} approach, i.e., a combination of analytical and 
numerical techniques. In the present case this would mean to 
introduce functions $\nu_{i}=e^{-f/\epsilon}\psi_{i}$, $i=1,2$, and solve 
numerically the system following from (\ref{eq:specprob2}) for the 
functions $\nu_{1,2}$ for a given solution $f$ to the eikonal 
equation. This will be the subject of further work.

\section*{Acknowledgement}
This work was partially supported by the PARI and FEDER programs in 
2016 and 2017, by the ANR-FWF project ANuI and by the Marie-Curie 
RISE network IPaDEGAN. K.M. was supported in part by the National 
Science Foundation under grant DMS-1733967. We are grateful to J. Sj\"ostrand for 
helpful discussions and hints.

\end{document}